\documentclass[titlepage,twoside,12pt]{article}
\usepackage{amssymb}
\usepackage{amsfonts}
\begin{document}

{\LARGE \bf Group Invariant Entanglements in \\ \\ Generalized Tensor Products} \\ \\


{\bf Elem\'{e}r E ~Rosinger} \\ \\
{\small \it Department of Mathematics \\ and Applied Mathematics} \\
{\small \it University of Pretoria} \\
{\small \it Pretoria} \\
{\small \it 0002 South Africa} \\
{\small \it eerosinger@hotmail.com} \\ \\

{\bf Abstract} \\

The group invariance of entanglement is obtained within a very general and simple setup of the latter, given by a recently introduced considerably
extended concept of tensor products. This general approach to entanglement - unlike the usual one given in the particular setup of tensor products of
vector spaces - turns out not to need any specific algebraic structure. The resulting advantage is that, entanglement being in fact defined by a
negation, its presence in a general setup increases the chances of its manifestations, thus also its availability as a resource. \\ \\ \\

{\bf 0. Preliminaries} \\

The interest in generalized tensor products, [3-5], is natural and also quite fundamental. Let us indeed divide physical systems in two classes : \\

1) Cartesian : those for which the state space of the composite is the Cartesian product of the state spaces of the components. Therefore obviously, Cartesian systems cannot exhibit entanglement. \\

2) Non-Cartesian : those for which the state space of the composite is not the Cartesian product of the state spaces of the components, and instead, it is some larger set. Therefore, such non-Cartesian systems do inevitably exhibit entanglement. \\

Regarding the above division, what is not yet understood clearly enough is that the usual quantum systems are but a particular instance of non-Cartesian systems. \\

The generalized tensor products provide a considerably large family of non-Cartesian ways of composition of systems. Indeed, so far, it was believed that only vector spaces could be composed by tensor products. Thus by far most of the non-quantum systems automatically fell out of the possibility of composition by tensor products, since their state spaces where not vector spaces. \\
However, generalized tensor products can now be defined for absolutely arbitrary sets, and all one needs is some very mild structures on them, far milder than any algebraic ones, and in particular, far milder than vector spaces. \\

And then the only problem is to find some really existing systems which compose at least in some of these many non-Cartesian ways. \\

A main interest in this regard comes, of course, from quantum computation. \\

Here however, we should better use the term non-Cartesian computation. Indeed, the electronic computers are of course Cartesian systems, and our great interest in quantum computers is not in the quanta themselves, but in their non-Cartesian composition, which among others, opens the possibility to such extraordinary computational resources as entanglement. \\

However, the immense trouble with the quantum type non-Cartesian computers is decoherence. \\

And then, the idea is to use non-Cartesian systems which do not decohere, yet have available entanglement. \\

Not to mention that, by studying lots of non-Cartesian systems other than the quantum ones, one may discover new and yet unknown resources, beyond entanglement. \\

In conclusion :

1) Non-Cartesian systems can be many more than quantum ones, and they all have entanglement. \\

2) Let us call Classical systems those which do NOT exhibit decoherence, or in which decoherence can easily be avoided. \\

And then, we are interested in : \\

Non-Cartesian Classical systems which, therefore, have entanglement and do not have decoherence. \\

The fact that there may be many many non-Cartesian systems is suggested by the immense generality and variety of generalized tensor products.  All that remains to do, therefore, is to identify as many as possible such really existing non-Cartesian Classical systems ... \\
And of course, also to identify important new features of them beyond entanglement ... \\ \\

{\bf 1. Group Actions on General Tensor Products} \\

We shall consider entanglement within the very general and rather simple underlying tensor product structure, briefly
presented for convenience in Appendix 2. \\

Let $X,~ Y$ be arbitrary nonvoid sets and let ${\cal A}$ be any nonvoid set of binary operations on $X$, while correspondingly, ${\cal B}$ is any nonvoid
set of binary operations on $Y$. Then, as seen in Appendix 2, one can define the tensor  product \\

(1.1)~~~ $ X \bigotimes_{{\cal A}, {\cal B}} Y = Z / \approx_{{\cal A}, {\cal B}} $ \\

with the canonical quotient embedding \\

(1.2)~~~ $ X \times Y \ni ( x, y ) \longmapsto x \bigotimes_{{\cal A}, {\cal B}} y \in X \bigotimes_{{\cal A}, {\cal B}} Y $ \\

and consequently $X \bigotimes_{{\cal A}, {\cal B}} Y$ will be the set of all elements \\

(1.3)~~~ $ x_1 \bigotimes_{{\cal A}, {\cal B}} y_1 ~\gamma~ x_2 \bigotimes_{{\cal A}, {\cal B}} y_2 ~\gamma~
                                  \ldots ~\gamma~ x_n \bigotimes_{{\cal A}, {\cal B}} y_n $ \\

with $n \geq 1$ and $x_i \in X,~ y_i \in Y$, for $1 \leq i \leq n$. \\

Let now $( G, . ),~ ( H, . )$ be two groups which act on $X$ and $Y$, respectively, according to \\

(1.4)~~~ $ G \times X \ni ( g , x ) \longmapsto g x \in X,~~~ H \times Y \ni ( h, y ) \longmapsto h y \in Y $ \\

It is important to note that, as a general property of such group actions, for every given $g \in G,~ h \in H$, the mappings \\

(1.5)~~~ $ X \ni x \longmapsto g x \in X,~~~ Y \ni y \longmapsto h y \in Y $ \\

are bijective. \\

Our aim is to define a group action $( G, . ) \times ( H, . )$ on, see (A2.2.4), the tensor product $X \bigotimes_{{\cal
A}, {\cal B}} Y$. \\

In this regard we note that, as a consequence of (1.4), we have a natural group action \\

(1.6)~~~ $ ( G \times H ) \times ( X \times Y ) \ni ( ( g, h ), ( x, y ) ) \longmapsto ( g x, h y ) \in X \times Y $ \\

thus we obtain a natural group action, see (A2.1.2) \\

(1.7)~~~ $ \begin{array}{l}
                ( G \times H ) \times Z \ni ( ( g, h ), ( ( x_1, y_1 ) ~\gamma~ ( x_2, y_2 ) ~\gamma~ \ldots ~\gamma~ ( x_n, y_n ) ) )
                          \longmapsto \\  \\
                   ~~~~~~~ \longmapsto ( g x_1, h y_1 ) ~\gamma~ ( g x_2, h y_2 ) ~\gamma~ \ldots ~\gamma~ ( g x_n, h y_n ) \in Z
            \end{array} $ \\

Now we shall assume that the families of binary operations ${\cal A},~ {\cal B}$ satisfy the following compatibility relations
with the respective group actions $( G, . ),~ ( H, . )$, namely \\

(1.8)~~~ $ \begin{array}{l}
                    \forall~~~ \alpha \in {\cal A},~~ \beta \in {\cal B} ~: \\ \\
                    \forall~~~ g \in G,~ h \in H,~ x, x\,' \in X,~ y, y\,' \in Y ~: \\ \\
                    ~~~~~ \alpha ( g x , g x\,' ) = g \alpha ( x, x\,' ),~~~ \beta ( h y, h y\,' ) = h \beta ( y, y\,' )
            \end{array} $ \\ \\

In this case, the equivalence relation, see (A2.2.1) - (A2.2.3), $\approx_{{\cal A}, {\cal B}}$ on $Z$ has, for $z, z\,'
\in Z$, and $g  \in G,~ h \in H$, the following property \\

(1.9)~~~ $ z ~\approx_{{\cal A}, {\cal B}}~ z\,' ~~~\Longleftrightarrow~~~ ( g, h ) z ~\approx_{{\cal A}, {\cal B}}~ ( g, h ) z\,' $ \\

Indeed, let us assume the left hand side of the above relation. Then in view of (A2.2.1) - (A2.2.3), we have three
possibilities, namely \\

(1.10)~~~ $ z = z\,' $ \\

(1.11)~~~ for some $\alpha \in {\cal A}$, the transformation (A2.2.1) applied to $z$ \\
          \hspace*{1.5cm} gives $z\,'$, or vice-versa \\

(1.12)~~~ for some $\beta \in {\cal B}$, the transformation (A2.2.2) applied to $z$ \\
          \hspace*{1.5cm} gives $z\,'$, or vice-versa \\

Now if (1.10) holds, then the right hand term of (1.9) is obviously valid. \\

Further, in the case of (1.11), we have \\

$~~~~~~ z = ( x_1, y_1 ) ~\gamma~ ( x\,'_1, y_1 ) ~\gamma~ ( x_2, y_2 ) ~\gamma~ \ldots ~\gamma~ ( x_n, y_n ) $ \\

$~~~~~~ z\,' = ( \alpha ( x_1, x\,'_1), y_1 ) ~\gamma~ ( x_2, y_2 ) ~\gamma~ \ldots , ( x_n, y_n ) $ \\

thus (1.7) gives  for $g \in G,~ h \in H$, \\

$~~~~~~ ( g, h ) z = ( g x_1, h y_1 ) ~\gamma~ ( g x\,'_1, h y_1 ) ~\gamma~ ( g x_2, h y_2 ) ~\gamma~
                                  \ldots ~\gamma~ ( g x_n, h y_n ) $ \\

$~~~~~~ ( g, h ) z\,' = ( g \alpha ( x_1, x\,'_1), h y_1 ) ~\gamma~ ( g x\,'_1, h y_1 ) ~\gamma~ ( g x_2, h y_2 )
                              ~\gamma~ \ldots ~\gamma~ ( g x_n, h y_n ) $ \\

with the second one, in view of (1.8), being \\

$~~~~~~ ( g, h ) z\,' = ( \alpha ( g x_1, g x\,'_1), h y_1 ) ~\gamma~ ( g x\,'_1, h y_1 ) ~\gamma~ ( g x_2, h y_2 )
                              ~\gamma~ \ldots ~\gamma~ ( g x_n, h y_n ) $ \\

Thus applying (1.11) once again, the right hand term of (1.9) again holds. \\

In case of (1.12), the argument is similar. \\

Conversely, let us assume the validity of the right hand term in (1.9). Then applying to it the group action $( g^{-1}, h^{-1} )$ and the above argument,
we obviously obtain the left hand term in (1.9). \\

Now the validity of (1.9) means that the group action (1.6) can naturally be extended, see (1.2), to a group action \\

(1.13)~~~ $ ( G \times H ) \times ( X \bigotimes_{{\cal A}, {\cal B}} Y )
                                \longrightarrow ( X \bigotimes_{{\cal A}, {\cal B}} Y ) $ \\

by \\

(1.14)~~~ $ ( g, h ) ( x \bigotimes_{{\cal A}, {\cal B}} y ) = ( g x ) \bigotimes_{{\cal A}, {\cal B}} ( h y ) $ \\

with the consequence \\

(1.15)~~~ $ \begin{array}{l}
             ( g, h ) ( x_1 \bigotimes_{{\cal A}, {\cal B}} y_1 ~\gamma~ x_2 \bigotimes_{{\cal A}, {\cal B}} y_2 ~\gamma~
             \ldots ~\gamma~ x_n \bigotimes_{{\cal A}, {\cal B}} y_n ) = \\ \\
             ~~~~ = ( g x_1 ) \bigotimes_{{\cal A}, {\cal B}} ( h y_1 ) ~\gamma~ ( g x_2 ) \bigotimes_{{\cal A}, {\cal B}}
             ( h y_2 ) ~\gamma~ \ldots ~\gamma~ ( g x_n ) \bigotimes_{{\cal A}, {\cal B}} ( h y_n )
            \end{array} $ \\

for $g \in G,~ h \in H$. \\ \\

{\bf 2. Group Actions on Generalized Entanglements} \\

We recall the general definition of entanglement within the above extended concept of tensor products. Namely, an element \\

(2.1)~~~ $ w \in X \bigotimes_{{\cal A}, {\cal B}} Y $ \\

is called {\it entangled}, if and only if it is {\it not} of the form \\

(2.2)~~~ $ w = x \bigotimes_{{\cal A}, {\cal B}} y $ \\

for some $x \in X$ and $y \in Y$. \\

{\bf Theorem 2.1.} \\

Given any group action, see (1.13) - (1.15) \\

(2.3)~~~ $ ( G \times H ) \times ( X \bigotimes_{{\cal A}, {\cal B}} Y )
                                \longrightarrow ( X \bigotimes_{{\cal A}, {\cal B}} Y ) $ \\

and element \\

(2.4)~~~ $ w = x \bigotimes_{{\cal A}, {\cal B}} y  X \bigotimes_{{\cal A}, {\cal B}} Y $ \\

then, for every $( g, h ) \in G \times H$, we have \\

(2.5)~~~ $ w $ is entangled ~~~$\Longleftrightarrow~~~ ( g, h ) w $ is entangled \\

therefore \\

(2.6)~~~ $ w $ is not entangled ~~~$\Longleftrightarrow~~~ ( g, h ) w $ is not entangled \\

{\bf Proof.} \\

Let us assume the situation in the left hand term of (2.6). Then (2.2) gives \\

$~~~~~~ w = x \bigotimes_{{\cal A}, {\cal B}} y $ \\

for some $x \in X$ and $y \in Y$, hence in view of (1.14), we have \\

$~~~~~~ ( g, h ) w = ( g x ) \bigotimes_{{\cal A}, {\cal B}}( h y ) $ \\

which according to (2.2), implies that  the right hand term of (2.6) is valid. \\

Obviously, the converse implication in (2.6) follows from the above line of argument, by applying the group action  $( g^{-1}, h^{-1} )$ to the right
hand term in it. \\

As for (2.5), it is but an immediate consequence of (2.6). \\

{\bf Remark 2.1.} \\

The relevance of Theorem 2.1. is in the following three facts : \\

1) The arbitrary group actions (1.13) - (1.15) \\

(2.7)~~~ $ ( G \times H ) \times ( X \bigotimes_{{\cal A}, {\cal B}} Y )
                                \longrightarrow ( X \bigotimes_{{\cal A}, {\cal B}} Y ) $ \\

on the general tensor products, see (1.1) \\

(2.8)~~~ $ X \bigotimes_{{\cal A}, {\cal B}} Y $ \\

are can be defined by only requiring the natural compatibility conditions, see (1.8) \\

(2.9)~~~  $ \begin{array}{l}
                    \forall~~~ \alpha \in {\cal A},~~ \beta \in {\cal B} ~: \\ \\
                    \forall~~~ g \in G,~ h \in H,~ x, x\,' \in X,~ y, y\,' \in Y ~: \\ \\
                    ~~~~~ \alpha ( g x , g x\,' ) = g \alpha ( x, x\,' ),~~~ \beta ( h y, h y\,' ) = h \beta ( y, y\,' )
            \end{array} $ \\ \\

2) The property of being, or alternatively, not being entangled within these general tensor products (2.8), is {\it invariant} under the general group
actions (2.9). \\

3) Since being entangled is defined by the {\it negation} of a relation, namely, (2.2), it follows that, quite likely, there are far more entangled
elements in a given tensor product, than there are non-entangled ones, such indeed being the case with the usual tensor products of vector spaces.
Consequently, when generalizing the concept of tensor products to the extent done in (2.8), while at the same time, keeping the corresponding
generalization of the concept of entanglement still as the negation of the same kind of relation, it is quite likely that the amount of entangled
elements in such generalized tensor products may significantly {\it increase}. And such an increase may be convenient, if we recall the extent to which
usual entanglement is a fundamental resource in quantum information theory. \\ \\

{\bf 3. Examples} \\

We present here, starting with section 3.2. below, several examples with the following two aims, namely, to show : \\

1) how much more general is the context in which tensor products and entanglement can be defined, than the usual one based on vector spaces, \\

2) the variety and novelty of both tensor products and entanglement even in the simplest cases which go beyond the usual concepts. \\

First however, several preliminary constructions are needed. \\

{\bf 3.1. Preliminaries} \\

{\bf 3.1.1. Binary Operations} \\

Given any nonvoid set $X$ and any binary operation $\alpha : X \times X \longrightarrow X$, we call a subset $A
\subseteq X$ to be $\alpha$-{\it stable}, if and only if \\

(3.1.1.1)~~~ $ x, y \in A ~~\Longrightarrow~~ \alpha ( x, y ) \in A $ \\

Obviously, $X$ is $\alpha$-stable, and the intersection of any family of $\alpha$-stable subsets is $\alpha$-stable.
Consequently, for every subset $A \subseteq X$, we can define the smallest $\alpha$-stable subset which contains it,
namely \\

(3.1.1.2)~~~ $ [ A ]_\alpha = \bigcap_{A \subseteq B,~ B ~\alpha-stable}~ B $ \\

Therefore, we can associate with $\alpha$ the mapping  $\psi_\alpha : {\cal P} ( X ) \longrightarrow {\cal P} ( X )$
defined by \\

(3.1.1.3)~~~ $ \psi_\alpha ( A ) = [ A ]_\alpha,~~~ A \subseteq X $ \\

In view of (3.1.2), we have \\

(3.1.1.4)~~~ $ \psi_\alpha ( \psi_\alpha ( A ) ) = \psi_\alpha ( A ) \supseteq A,~~~ A \subseteq X $ \\

since as mentioned, $[ A ]_\alpha$ is $\alpha$-stable, and obviously $[ A ]_\alpha \subseteq [ A ]_\alpha$. \\

A particular case of the above is the following. Let $( S,\ast )$ be a semigroup with the neutral element $e$. Then
$[ \{ e \} ]_\ast = \{ e \}$, while for $a \in S,~ a \neq e$, we have $[ \{ a \} ]_\ast = \{ a, a \ast a, a \ast a \ast a,
\dots \}$. \\

For instance, if $( S, \ast ) = ( \mathbb{N}, + )$, then $[ \{ 0 \} ]_+ = \{ 0 \}$, while $[ \{ 1 \} ]_+ = \mathbb{N}
\setminus \{ 0 \} = \mathbb{N}_1$. \\

We note that the mapping $\psi_\alpha : {\cal P} ( X ) \longrightarrow {\cal P} ( X )$ satisfies three of the four
Kuratowski closure axioms, except for $\psi_\alpha ( A \cup B ) = \psi_\alpha ( A ) \cup \psi_\alpha ( B )$, with $A, B
\subseteq X$. Indeed, this condition is obviously not satisfied, as can be seen in the following simple example, when $X =
\mathbb{R}^2$ and $\alpha$ is the usual addition $+$. Let $A = \mathbb{}R \times \{ 0 \},~ B = \{ 0 \} \times \mathbb{R}
\subsetneqq \mathbb{R}^2$. Then $\psi_\alpha ( A ) = A,~ \psi_\alpha ( B ) = B$, while $\psi_\alpha ( A \cup B ) =
\mathbb{R}^2$. \\

We shall denote by \\

(3.1.1.5)~~~ $ {\cal B}_X $ \\

the set of all binary operations on $X$. Obviously \\

(3.1.1.6)~~~ $ car\, X = n ~~\Longrightarrow~~ car\, {\cal B}_X = n^{(n^2)} $ \\ \\

{\bf 3.1.2. Generators} \\

Given any nonvoid set $X$, a {\it generator} on $X$ is any mapping $\psi : {\cal P} ( X ) \longrightarrow {\cal P} ( X )$,
such that \\

(3.1.2.1)~~~ $ A \subseteq \psi ( A ),~~~ A \subseteq X $ \\

(3.1.2.2)~~~ $ \psi ( A ) \subseteq \psi ( A\,' ),~~~ A \subseteq A\,' \subseteq X $ \\

It follows that every $\psi_\alpha$ associated in (3.1.1.3) to a binary operation $\alpha$ on $X$ is a generator. \\ \\

{\bf 3.1.3. From Generators to Binary Operations} \\

Let us consider the inverse problem, namely, to associate binary operations to given generators. \\

Let $\psi$ be a generator on a set $X$. A binary operation $\alpha$ on $X$ is {\it compatible} with $\psi$, if and only
if \\

(3.1.3.1)~~~ $ \psi ( A ) $ is $\alpha$-stable,~~~ $ A \subseteq X $ \\

In view of (3.1.1.2), this condition is equivalent with \\

(3.1.3.2)~~~ $ [ A ]_\alpha \subseteq \psi ( A ),~~~ A \subseteq X $ \\

which is further equivalent with \\

(3.1.3.3)~~~ $ \psi_\alpha ( A ) \subseteq \psi ( A ),~~~ A \subseteq X $ \\

Let us now denote by \\

(3.1.3.4)~~~ $ {\cal B}_\psi $ \\

the set of all binary operations $\alpha$ on $X$ which are compatible with $\psi$. \\ \\

{\bf 3.1.4. Open Problems} \\

1) Given  a generator $\psi$ on a set $X$, find the binary operations $\alpha$ on $X$ compatible with $\psi$ and with the
largest $\psi_\alpha$ in the sense of (4.3). \\

2) Characterize those generators $\psi$ on a set $X$ which coincide with such largest $\psi_\alpha $. \\ \\

{\bf 3.1.5. Special Classes of Generators} \\ \\

{\bf 3.1.5.1. The Identity Generator} \\

We note that the identity mapping $id_{{\cal P} ( X )} : {\cal P} ( X ) \ni A \longmapsto A \in {\cal P} ( X )$ is a generator on $X$. Therefore, let us
start by presenting binary operations on $X$ which belong to ${\cal B}_{id_{{\cal P} ( X )} }$. \\

Further, let us denote by $\lambda_X,~ \rho_X$ the binary operations on $X$ defined by \\

(3.1.5.1.1)~~~ $ \lambda_X ( x, x\,' ) = x,~~~ \rho_X ( x, x\,' ) = x\,',~~~ x, x\,' \in X $ \\

then it is obvious that \\

(3.1.5.1.2)~~~ $ \psi_{\lambda_X} = \psi_{\rho_X} = id_{{\cal P} ( X )} $ \\

Furthermore, let ${\cal Q} \subseteq X \times X$ and define $\alpha_{\cal Q} : X \times X \longrightarrow X$ by \\

(3.1.5.1.3)~~~ $ \alpha_{\cal Q} ( x, x\,' ) = \begin{array}{|l}
                                                      x ~\mbox{if}~ ( x, x\,' ) \in {\cal Q} \\ \\
                                                      x\,' ~\mbox{if}~ ( x, x\,' ) \notin {\cal Q}
                                          \end{array} $ \\

then again, we obtain \\

(3.1.5.1.4)~~~ $ \psi_{\alpha_{\cal Q}} = id_{{\cal P} ( X )} $ \\

Let now $\leq$ be a total order on $X$ and denote by $\min_X,~ \max_X$ the binary operations on $X$ defined by \\

(3.1.5.1.5)~~~ $ \min_X ( x, x\,' ) = x \wedge x\,',~~~ \max_X ( x, x\,' ) = x \vee x\,',~~~ x, x\,' \in X $ \\

Obviously once more we have \\

(3.1.5.1.6)~~~ $ \psi_{\min_X} = \psi_{\max_X} = id_{{\cal P} ( X )} $ \\

{\bf Remark 3.1.5.1.1.} \\

We note that, see (3.1.5.1.3) \\

(3.1.5.1.7)~~~ $\lambda_X = \alpha_{X \times X},~~~ \rho_X = \alpha_\phi $ \\

(3.1.5.1.8)~~~ $ \min_X = \alpha_{X^2_+},~~~ \max_X = \alpha_{X^2_-} $ \\

where $X^2_+ = \{ ( x, x\,' ) \in X \times X ~|~ x \leq x\,' \}$, while $X^2_- = \{ ( x, x\,' ) \in X \times X ~|~ x\,' \leq x \}$. \\

We also note that both binary operations $\lambda_X$ and $\rho_X$ are {\it associative}, while in case $X$ has at least two elements, they are {\it not}
commutative. Similarly, for a total order $( X, \leq )$, both binary operations $\min_X$ and $\max_X$ are {\it associative}, and when $X$ has at least two
elements, they are {\it not} commutative. \\

{\bf Theorem 3.1.5.1.1.} \\

Let $\alpha$ be a binary operation on $X$. Then the following three properties are equivalent \\

(3.1.5.1.9)~~~ $ \psi_\alpha = id_{{\cal P} ( X )} $ \\

(3.1.5.1.10)~~~ $ \forall~ x, x\,' \in X ~: ~~ \alpha ( x, x\,' ) \in \{ x, x\,'\} $ \\

(3.1.5.1.11)~~~ $ \exists~ {\cal Q} \subseteq X \times X ~: ~~ \alpha = \alpha_{\cal Q} $ \\

{\bf Proof.} \\

In view of (3.1.5.1.4), we have (3.1.5.1.11) $\Longrightarrow$ (3.1.5.1.9). \\
Assume $\psi_\alpha = id_{{\cal P} ( X )}$, then (3.1.3) yields $[ A ]_\alpha = A$, for $A \subseteq X$. Given now $A = \{ x, x\,' \} \subseteq X$, it
follows that $\alpha ( x, x\,' ) \in \{ x, x\,'\}$. Thus (3.1.5.1.9) $\Longrightarrow$ (3.1.5.1.10). \\
We show now that (3.1.5.1.10) $\Longrightarrow$ (3.1.5.1.11). Let \\

$~~~~~~ {\cal Q} = \{~ ( x, x\,' ) \in X \times X ~|~ \alpha ( x, x\,' ) = x ~\} $ \\

then clearly $\alpha = \alpha_{\cal Q}$.

\hfill $\Box$ \\

As seen next, there are plenty of binary operations $\alpha$ on a set $X$, such that $\psi_\alpha = id_{{\cal P} ( X )}$. \\

{\bf Theorem 3.1.5.1.2.} \\

Given ${\cal Q},~ {\cal Q}\,' \subseteq X \times X$, then \\

(3.1.5.1.12)~~~ $ \alpha_{\cal Q} = \alpha_{{\cal Q}\,'}
                      ~~~\Longleftrightarrow~~~ {\cal Q} \setminus \Delta^2 X = {\cal Q}\,' \setminus \Delta^2 X $ \\

Consequently \\

(3.1.5.1.13)~~~ $ car\, X = n ~~\Longrightarrow~~ car\, {\cal B}_{id_{{\cal P} ( X )}} = 2^{( n^2 - n)} $ \\

{\bf Proof.} \\

The equivalence in (3.1.5.1.12) is immediate. Further, one notes that (3.1.5.1.12) $\Longrightarrow$ (3.1.5.1.13).

\hfill $\Box$ \\

Let us find now all associative binary operations in ${\cal B}_{id_{{\cal P} ( X )} }$. The answer is presented in \\

{\bf Theorem 3.1.5.1.3.} \\

Given ${\cal Q} \subseteq X \times X$, then \\

(3.1.5.1.14)~~~ $ \alpha_{\cal Q} $ is associative
                     ~~~$\Longleftrightarrow~~~ {\cal Q} \circ {\cal Q} \subseteq {\cal Q} $ \\

{\bf Proof.} \\

Let us start with the implication $\Longrightarrow$. Given $x, x\,', x\,'' \in X$, assume that $( x, x\,' ),
( x\,', x\,'' ) \in {\cal Q}$. Then $\alpha_{\cal Q} ( \alpha_{\cal Q} ( x, x\,' ), x\,'' ) =
\alpha_{\cal Q} ( x, x\,'' )$, while $\alpha_{\cal Q} ( x, \alpha_{\cal Q} ( x\,', x\,'' ) ) =
\alpha_{\cal Q} ( x, x\,' )= x$. Thus the associativity of $\alpha_{\cal Q}$ implies $\alpha_{\cal Q} ( x, x\,'' ) = x$,
which means that $( x, x\,'' ) \in {\cal Q}$. \\

For the converse implication $\Longleftarrow$, let $x, x\,', x\,'' \in X$. Then we have the following four possible
cases : \\

1) $( x, x\,' ), ( x\,', x\,'' ) \in {\cal Q}$ \\

This means that $( x, x\,'' ) \in {\cal Q}$, thus $\alpha_{\cal Q} ( \alpha_{\cal Q} ( x, x\,' ), x\,'' ) =
\alpha_{\cal Q} ( x, x\,'' ) = x$. \\
On the other hand, $\alpha_{\cal Q} ( x, \alpha_{\cal Q} ( x\,', x\,'' ) ) =  \alpha_{\cal Q} ( x, x\,' ) = x$. Hence
associativity holds. \\

2) $( x, x\,' ) \in {\cal Q},~ ( x\,', x\,'' ) \notin {\cal Q}$ \\

Then $\alpha_{\cal Q} ( \alpha_{\cal Q} ( x, x\,' ), x\,'' ) = \alpha_{\cal Q} ( x, x\,'' )$, while
$\alpha_{\cal Q} ( x, \alpha_{\cal Q} ( x\,', x\,'' ) ) = \alpha_{\cal Q} ( x, x\,'' )$, thus associativity holds. \\

3) $( x, x\,' ), ( x\,', x\,'' ) \notin {\cal Q}$ \\

Then $\alpha_{\cal Q} ( \alpha_{\cal Q} ( x, x\,' ), x\,'' ) = \alpha_{\cal Q} ( x\,', x\,'' ) = x\,''$, while
$\alpha_{\cal Q} ( x, \alpha_{\cal Q} ( x\,', x\,'' ) ) = \alpha_{\cal Q} ( x, x\,'' )$. Thus associativity holds, if and only if \\

(3.1.5.1.15)~~~ $ \alpha_{\cal Q} ( x, x\,'' ) = x\,'' $ \\

which in view of (3.1.5.1.3) is equivalent with \\

(3.1.5.1.16)~~~ $ ( x, x\,'' ) \notin {\cal Q} $ \\

However, we note that, for $u, v \in X,~ u \neq v$, we have \\

(3.1.5.1.17)~~~ $ ( u, v ) \notin {\cal Q} ~~~\Longleftrightarrow~~~ ( v, u ) \in {\cal Q} $ \\

and in view of (3.1.5.1.12), we can assume that \\

(3.1.5.1.18)~~~ $ \Delta^2 X \subseteq {\cal Q} $ \\

But then, the assumption $( x, x\,' ), ( x\,', x\,'' ) \notin {\cal Q}$ implies $x \neq x\,',~ x\,' \neq x\,''$, hence (3.1.5.1.17) gives $( x\,'', x\,' ),
( x\,', x ) \in {\cal Q}$, which in view of the assumed right hand in (3.1.5.1.14), means \\

(3.1.5.1.19)~~~ $ ( x\,'', x ) \in {\cal Q} $ \\

Here, there are another two subcases, namely \\

3.1) $x\,'' \neq x$ \\

And this, together with (3.1.5.1.19), (3.1.5.1.17), gives (3.1.5.1.16), thus (3.1.5.1.15), which yields associativity. \\

3.2) $x\,'' = x$ \\

In which case (3.1.5.1.15) follows trivially, thus also associativity. \\

4) $( x, x\,' ) \notin {\cal Q},~ ( x\,', x\,'' ) \in {\cal Q}$ \\

Then $\alpha_{\cal Q} ( \alpha_{\cal Q} ( x, x\,' ), x\,'' ) = \alpha_{\cal Q} ( x\,', x\,'' ) = x\,'$, while
$\alpha_{\cal Q} ( x, \alpha_{\cal Q} ( x\,', x\,'' ) ) = \alpha_{\cal Q} ( x, x\,' )= x\,'$, thus associativity holds.

\hfill $\Box$ \\

As for the commutative binary operations in ${\cal B}_{id_{{\cal P} ( X )} }$, we have the obvious result in \\

{\bf Theorem 3.1.5.1.4.} \\

Given ${\cal Q} \subseteq X \times X$, then \\

(3.1.5.1.20)~~~ $ \alpha_{\cal Q} $ is commutative ~~~$\Longleftrightarrow~~~ {\cal Q} = \Delta^2 X $ \\

{\bf Remark 3.1.5.1.2.} \\

In the particular case of the identity generators $id_{{\cal P} ( X )}$ on arbitrary sets $X$, the two problems in section 5 above obtain simple
solutions. Namely, every binary operation $\alpha \in {\cal B}_{id_{{\cal P} ( X )}}$ is a solution of both mentioned problems. \\

Indeed, according to (3.1.5.1.11), each such binary operation is of t he form $\alpha = \alpha_{\cal Q}$, for some ${\cal Q} \subseteq X \times X$. And then
(3.1.5.1.9) concludes the proof. \\ \\

{\bf 3.1.5.2. A Few Steps beyond the Identity Generator} \\

Generators $\psi$ on $X$ which are more complex than the identity generator can obviously be obtained in many ways. In the sequel, we shall consider
several such classes. \\

A first possible class of generators $\psi$ on $X$ which are more complex than the identity generator is given by those for which \\

(3.1.5.2.1)~~~ $ \exists~~~ \alpha \in {\cal B}_\psi ~:~~~ \exists~~~ x, x\,' \in X ~:~~~ \alpha ( x, x\,' ) \notin \{ x, x\,' \} $ \\

A simple example can be obtained as follows. Let $a, b, c \in X$ be three different elements, and let $\alpha : X \times X \longrightarrow X$ be given
by \\

(3.1.5.2.2)~~~ $ \alpha ( x, x\,' ) = \begin{array}{|l}
                                                      x ~\mbox{if}~ ( x, x\,' ) \neq ( a, b ) \\ \\
                                                      c ~\mbox{if}~ ( x, x\,' ) = ( a, b )
                                          \end{array} $ \\

then $\alpha$ is neither associative, nor commutative. \\ \\

{\bf 3.1.6. Open Problems} \\

1) Characterize the generators $\psi$ on $X$ for which (3.1.5.2.1) holds with $\alpha$ associative. \\

2) Characterize the generators $\psi$ on $X$ for which (3.1.5.2.2) holds with $\alpha$ commutative. \\ \\

{\bf 3.2. A First Example : Tensor Products of Totally Ordered \\
          \hspace*{1cm }Sets, and Their Entanglement} \\

Let $( X, \leq )$ and $( Y, \leq)$ be two totally ordered sets, and let us consider on them the respective binary operations, see (3.1.5.1.5), $\min_X$
and $\min_Y$. We shall now construct the tensor product, see (A2.1.8) \\

(3.2.1)~~~ $ X \bigotimes_{\min_X, \min_Y} Y = Z / \approx_{\min_X, \min_Y} $ \\

For that purpose, we have to particularize the conditions (A2.1.6), (A2.1.7) which are involved in the definition of the equivalence relation
$\approx_{\min_X, \min_Y}$. In this respect, obviously (A2.1.6) takes the from \\

(3.2.2)~~~  replace $( x_1, y_1 ) ~\gamma~ ( x_2, y_2 ) ~\gamma~ \ldots ~\gamma~ ( x_n, y_n )$ \\
            \hspace*{1.85cm} with $( x_1, y_1 ) ~\gamma~ ( x\,'_1, y_1 ) ~\gamma~ ( x_2, y_2 ) ~\gamma~ \ldots ~\gamma~ ( x_n, y_n )$, \\
            \hspace*{1.9cm}or vice-versa, where $x\,'_1 \in X,~ x\,'_1 \geq x_1$ \\

while (A2.1.7) becomes \\

(3.2.3)~~~  replace $( x_1, y_1 ) ~\gamma~ ( x_2, y_2 ) ~\gamma~ \ldots ~\gamma~ ( x_n, y_n )$ \\
            \hspace*{1.85cm} with $( x_1, y_1 ) ~\gamma~ ( x_1, y\,'_1 ) ~\gamma~ ( x_2, y_2 ) ~\gamma~ \ldots ~\gamma~ ( x_n, y_n )$, \\
            \hspace*{1.9cm}or vice-versa, where $y\,'_1 \in Y,~ y\,'_1 \geq y_1$ \\

It follows in particular that, for $x \in X,~ y \in Y$, we have \\

(3.2.4)~~~ $ \begin{array}{l}
                  x \bigotimes_{\min_X, \min_Y} y ~=~ ( x \bigotimes_{\min_X, \min_Y} y ) ~\gamma~ ( x\,' \bigotimes_{\min_X, \min_Y} y ) ~=~ \\ \\
                  ~=~ ( x \bigotimes_{\min_X, \min_Y} y ) ~\gamma~ ( x \bigotimes_{\min_X, \min_Y} y\,' )
              \end{array} $ \\ \\

for $x\,' \in X,~ x\,' \geq x,~ y\,' \in Y,~ y\,' \geq y$. Indeed, (3.2.2) gives \\

$~~~~~~ ( x, y ) ~\approx_{\min_X, \min_Y}~ ( x, y ) ~\gamma~ ( x\,', y ) $ \\

while (3.2.3) implies \\

$~~~~~~ ( x, y ) ~\approx_{\min_X, \min_Y}~ ( x, y ) ~\gamma~ ( x, y\,' ) $ \\

which proves (3.2.4). \\

We note that (3.2.4) is a powerful {\it simplification} rule for elements of the tensor product $X \bigotimes_{\min_X, \min_Y} Y$, which according to
(A2.1.10), are of the form \\

(3.2.5)~~~ $ x_1 \bigotimes_{\min_X, \min_Y} y_1 ~\gamma~ x_2 \bigotimes_{\min_X, \min_Y} y_2 ~\gamma~
                                                             \ldots ~\gamma~ x_n \bigotimes_{\min_X, \min_Y} y_n $ \\

The interpretation of (3.2.4) in terms of certain {\it paths} in  $( X, \leq ) \times ( Y, \leq )$ further clarifies its meaning. Indeed, let us
consider on $X \times Y$ the {\it reflexive, antisymmetric} binary relation $\dashv$, defined by \\

(3.2.6)~~~ $ ( x, y ) \dashv ( x \,', y\,' ) ~~~\Longleftrightarrow~~~
                             \left ( ~ \begin{array}{l}
                                            ~~1)~~ x = x\,',~~~ y \leq y\,',~~ \mbox{or} \\ \\
                                            ~~2)~~ x \leq x\,',~~~ y = y\,'
                                        \end{array} ~ \right ) $ \\ \\

which in general is {\it not} transitive. Then a finite sequence in $X \times Y$ \\

(3.2.7)~~~ $ ( x_1, y_1 ), ( x_2, y_2 ), ( x_3, y_3 ), \ldots , ( x_n, y_n ) $ \\

is called {\it path-free} in $( X, \leq ) \times ( Y, \leq )$, if and only if none of the relations holds \\

(3.2.8)~~~ $ ( x_i, y_i ) \dashv ( x_j, y_j ) $ \\

where $1 \leq i, j \leq n,~ i \neq j$. It is, therefore, convenient to introduce on $X \times Y$ the related binary relation $\bowtie$ defined by \\

(3.2.9)~~~ $ ( x, y ) \bowtie ( x \,', y\,' ) ~~~\Longleftrightarrow~~~
                              \left (~ \begin{array}{l}
                                          \mbox{neither}~~ ( x, y ) \dashv ( x\,', y\,' ) \\ \\
                                          \mbox{nor}~~ ( x\,', y\,' ) \dashv ( x, y )
                                        \end{array} ~ \right ) $ \\

It follows that the path-free condition (3.2.8) can be written in the equivalent form \\

(3.2.10)~~~ $ ( x_i, y_i ) \bowtie ( x_j, y_j ),~~~ 1 \leq i < j \leq n $ \\

Now, in view of (3.2.4), we have for $x, x\,' \in X,~ y, y\,' \in Y$, the implication \\

(3.2.11)~~~ $ ( x, y ) \dashv ( x \,', y\,' ) ~~~\Longrightarrow~~~
                                     \left ( ~ \begin{array}{l}
                                                   x \bigotimes_{\min_X, \min_Y} y ~=~ \\ \\
                                                   ~=~ ( x \bigotimes_{\min_X, \min_Y} y ) ~\gamma~ ( x\,' \bigotimes_{\min_X, \min_Y} y\,' )
                                                \end{array} ~ \right ) $ \\

and therefore \\

(3.2.12)~~~ $ \begin{array}{l}
               \left ( ~ \begin{array}{l}
                                x \bigotimes_{\min_X, \min_Y} y ~\neq~ \\ \\
                                ~=~ ( x \bigotimes_{\min_X, \min_Y} y ) ~\gamma~ ( x\,' \bigotimes_{\min_X, \min_Y} y\,' )
                          \end{array} ~ \right ) ~~~\Longrightarrow~~~ \\ \\

                   ~~~~~~~~~~~~ \Longrightarrow~~~ ( x, y ) \bowtie ( x \,', y\,' )
              \end{array} $ \\

And then (3.2.5) results in \\

{\bf Theorem 3.2.1.} \\

The elements of the tensor product  $X \bigotimes_{\min_X, \min_Y} Y$ are of the form \\

(3.2.13)~~~ $ x_1 \bigotimes_{\min_X, \min_Y} y_1 ~\gamma~ x_2 \bigotimes_{\min_X, \min_Y} y_2 ~\gamma~
                                   \ldots ~\gamma~ x_n \bigotimes_{\min_X, \min_Y} y_n $ \\

where \\

(3.2.14)~~~ $ ( x_1, y_1 ), ( x_2, y_2 ), ( x_3, y_3 ), \ldots , ( x_n, y_n ) $ \\

are path-free in $( X, \leq ) \times ( Y, \leq )$. \\

{\bf Corollary 3.2.1.} \\

If an element $w \in X \bigotimes_{\min_X, \min_Y} Y$ is entangled, then it is of form (3.2.13), and the respective path-free (3.2.14) is of length
at least two. \\

{\bf Remark 3.2.1.} \\

Obviously, the tensor products \\

$~~~~~~ X \bigotimes_{\min_X, \max_Y} Y,~ X \bigotimes_{\max_X, \min_Y} Y $ and $ X \bigotimes_{\max_X, \max_Y} Y $ \\

can be obtained in a similar manner, simply by replacing $( X, \leq )$ with $( X, \geq )$, and $( Y, \leq )$ with
$( Y, \geq )$. \\ \\

{\bf 3.3. A Second Example} \\

We can consider the more general situation when $( X, \leq )$ and $( Y, \leq)$ are two partially ordered sets which are
lattices. In this case we can still define the binary operations $\min_X : X \times X \longrightarrow X$ and $\min_Y :
Y \times Y \longrightarrow Y$ in the usual manner. \\

Let us now now construct the tensor product, see (A2.1.8) \\

(3.3.1)~~~ $ X \bigotimes_{\min_X, \min_Y} Y = Z / \approx_{\min_X, \min_Y} $ \\

For that purpose, we have to particularize the conditions (A2.1.6), (A2.1.7) which are involved in the definition of the equivalence relation
$\approx_{\min_X, \min_Y}$. In this respect, obviously (A2.1.6) takes the from \\

(3.3.2)~~~  replace $( x_1, y_1 ) ~\gamma~ ( x_2, y_2 ) ~\gamma~ \ldots ~\gamma~ ( x_n, y_n )$ \\
            \hspace*{1.85cm} with $( x'_1, y_1 ) ~\gamma~ ( x\,''_1, y_1 ) ~\gamma~ ( x_2, y_2 ) ~\gamma~
            \ldots ~\gamma~ ( x_n, y_n )$, \\
            \hspace*{1.9cm}or vice-versa, where $x\,'_1, x\,''_1 \in X,~ x = \min_X ( x\,'_1, x\,''_1 )$ \\

while (A2.1.7) becomes \\

(3.3.3)~~~  replace $( x_1, y_1 ) ~\gamma~ ( x_2, y_2 ) ~\gamma~ \ldots ~\gamma~ ( x_n, y_n )$ \\
            \hspace*{1.85cm} with $( x_1, y\,'_1 ) ~\gamma~ ( x_1, y\,''_1 ) ~\gamma~ ( x_2, y_2 ) ~\gamma~
            \ldots ~\gamma~ ( x_n, y_n )$, \\
            \hspace*{1.9cm}or vice-versa, where $y\,'_1, y\,''_1 \in Y,~ y = \min_Y ( y\,'_1, y\,''_1 )$ \\

It follows in particular that, for $x \in X,~ y \in Y$, we have \\

(3.3.4)~~~ $ \begin{array}{l}
                  x \bigotimes_{\min_X, \min_Y} y ~=~ ( x\,' \bigotimes_{\min_X, \min_Y} y )
                                     ~\gamma~ ( x\,'' \bigotimes_{\min_X, \min_Y} y ) ~=~ \\ \\
                  ~=~ ( x \bigotimes_{\min_X, \min_Y} y\,' ) ~\gamma~ ( x \bigotimes_{\min_X, \min_Y} y\,'' )
              \end{array} $ \\ \\

for $x\,', x\,'' \in X,~ x = \min_X ( x\,', x\,'' )$ and $y\,', y\,'' \in Y,~ y = \min_Y ( y\,', y\,'' )$. \\

Let us see the effect on the elements \\

(3.3.5)~~~ $ x_1 \bigotimes_{\min_X, \min_Y} y_1 ~\gamma~ x_2 \bigotimes_{\min_X, \min_Y} y_2 ~\gamma~
                                   \ldots ~\gamma~ x_n \bigotimes_{\min_X, \min_Y} y_n $ \\

which constitute the tensor product $X \bigotimes_{\min_X, \min_Y} Y$, of the simplifications introduced by
(3.3.4). For that purpose, the following notions are useful. Given $( x, y ) \in X \times Y$, we denote \\

(3.3.6)~~~ $ V_X ( x, y ) = \{~ ( ( x\,', y ), ( x\,'', y ) ) ~~|~~
                                                x\,', x\,'' \in X,~ x = \min_X ( x\,', x\,'' ) ~\} $ \\

and call it the $X$-{\it wedge} at $( x, y)$, while the $Y$-{\it wedge} at $( x, y)$ is given by \\

(3.3.7)~~~ $ V_Y ( x, y ) = \{~ ( ( x, y\,' ), ( x, y\,'' ) ) ~~|~~
                                                y\,', y\,'' \in Y,~ y = \min_Y ( y\,', y\,'' ) ~\} $ \\

Further, three pairs $( x_1, y_1 ), ( x_2, y_2 ), ( x_3, y_3 ) \in X \times Y$ are called {\it wedge-related}, if  and
only if, in one of their permutation $i, j, k$, we have \\

(3.3.8)~~~ $ ( ( x_i, y_i ), ( x_j, y_j ) ) \in V_X ( x_k, y_k ) \cup V_Y ( x_k, y_k ) $ \\

Now, a finite sequence in $X \times Y$ \\

(3.3.9)~~~ $ ( x_1, y_1 ), ( x_2, y_2 ), ( x_3, y_3 ), \ldots , ( x_n, y_n ) $ \\

is called {\it wedge-free} in $( X, \leq ) \times ( Y, \leq )$, if and only if no three of them are wedge-related. \\

It follows that (3.3.5) results in \\

{\bf Theorem 3.3.1.} \\

The elements of the tensor product  $X \bigotimes_{\min_X, \min_Y} Y$ are of the form \\

(3.3.10)~~~ $ x_1 \bigotimes_{\min_X, \min_Y} y_1 ~\gamma~ x_2 \bigotimes_{\min_X, \min_Y} y_2 ~\gamma~
                                   \ldots ~\gamma~ x_n \bigotimes_{\min_X, \min_Y} y_n $ \\

where \\

(3.3.11)~~~ $ ( x_1, y_1 ), ( x_2, y_2 ), ( x_3, y_3 ), \ldots , ( x_n, y_n ) $ \\

are wedge-free in $( X, \leq ) \times ( Y, \leq )$. \\

{\bf Corollary 3.3.1.} \\

If an element $w \in X \bigotimes_{\min_X, \min_Y} Y$ is entangled, then it is of form (3.3.10), and the respective
wedge-free (3.3.11) is of length at least two. \\

{\bf Remark 3.3.1.} \\

Obviously, the tensor products \\

$~~~~~~ X \bigotimes_{\min_X, \max_Y} Y,~ X \bigotimes_{\max_X, \min_Y} Y $ and $ X \bigotimes_{\max_X, \max_Y} Y $ \\

can be obtained in a similar manner, simply by replacing $( X, \leq )$ with $( X, \geq )$, and $( Y, \leq )$ with
$( Y, \geq )$. \\ \\

{\bf 3.4. A Third Example} \\

Let $X$ and $Y$ be two nonvoid sets and let us consider on them the binary operations, see (3.1.5.1.1), $\lambda_X,~ \lambda_Y$, respectively. Let us
now construct the tensor product, see (A2.1.8) \\

(3.4.1)~~~ $ X \bigotimes_{\lambda_X, \lambda_Y} Y = Z / \approx_{\lambda_X, \lambda_Y} $ \\

For that purpose, we have to particularize the conditions (A2.1.6), (A2.1.7) which are involved in the definition of the equivalence relation
$\approx_{\lambda_X, \lambda_Y}$. In this respect, obviously (A2.1.6) takes the from \\

(3.4.2)~~~  replace $( x_1, y_1 ) ~\gamma~ ( x_2, y_2 ) ~\gamma~ \ldots ~\gamma~ ( x_n, y_n )$ \\
            \hspace*{1.85cm} with $( x_1, y_1 ) ~\gamma~ ( x\,'_1, y_1 ) ~\gamma~ ( x_2, y_2 ) ~\gamma~ \ldots ~\gamma~ ( x_n, y_n )$, \\
            \hspace*{1.9cm}or vice-versa \\

while (A2.1.7) becomes \\

(3.4.3)~~~  replace $( x_1, y_1 ) ~\gamma~ ( x_2, y_2 ) ~\gamma~ \ldots ~\gamma~ ( x_n, y_n )$ \\
            \hspace*{1.85cm} with $( x_1, y_1 ) ~\gamma~ ( x_1, y\,'_1 ) ~\gamma~ ( x_2, y_2 ) ~\gamma~ \ldots ~\gamma~ ( x_n, y_n )$, \\
            \hspace*{1.9cm}or vice-versa \\

In particular, it follows that, for $x \in X,~ y \in Y$, we have \\

(3.4.4)~~~ $ \begin{array}{l}
                  x \bigotimes_{\lambda_X, \lambda_Y} y ~=~ ( x \bigotimes_{\lambda_X, \lambda_Y} y ) ~\gamma~
                               ( x\,' \bigotimes_{\lambda_X, \lambda_Y} y ) ~=~ \\ \\
                  ~=~ ( x \bigotimes_{\lambda_X, \lambda_Y} y ) ~\gamma~ ( x \bigotimes_{\lambda_X, \lambda_Y} y\,' )
              \end{array} $ \\ \\

for $x\,' \in X,~ y\,' \in Y$. Indeed, (3.4.2) gives \\

$~~~~~~ ( x, y ) ~\approx_{\lambda_X, \lambda_Y}~ ( x, y ) ~\gamma~ ( x\,', y ) $ \\

while (3.4.3) implies \\

$~~~~~~ ( x, y ) ~\approx_{\lambda_X, \lambda_Y}~ ( x, y ) ~\gamma~ ( x, y\,' ) $ \\

which gives the proof of (3.4.4). \\

Here we note that the above simplification rule (3.4.4) regarding the terms of the tensor product $X \bigotimes_{\lambda_X, \lambda_Y} Y$ is
considerably more powerful than that in (3.2.4). Consequently, the tensor product $X \bigotimes_{\lambda_X, \lambda_Y} Y$ will have {\it fewer} elements
than the tensor product $X \bigotimes_{\min_X, \min_Y} Y$. \\

Let us make more clear the respective situation. Any given finite sequence in $X \times Y$ \\

(3.4.5)~~~ $ ( x_1, y_1 ), ( x_2, y_2 ), ( x_3, y_3 ), \ldots , ( x_n, y_n ) $ \\

is called {\it repetition free}, if and only if \\

(3.4.6)~~~ $ x_i \neq x_j,~~ y_i \neq y_j,~~~ 1 \leq i < j \leq n $ \\

Now, in view of (A2.1.10), we obtain \\

{\bf Theorem 3.4.1.} \\

The elements of the tensor product  $X \bigotimes_{\lambda_X, \lambda_Y} Y$ are of the form \\

(3.4.7)~~~ $ x_1 \bigotimes_{\lambda_X, \lambda_Y} y_1 ~\gamma~ x_2 \bigotimes_{\lambda_X, \lambda_Y} y_2 ~\gamma~
                                   \ldots ~\gamma~ x_n \bigotimes_{\lambda_X, \lambda_Y} y_n $ \\

where \\

(3.4.8)~~~ $ ( x_1, y_1 ), ( x_2, y_2 ), ( x_3, y_3 ), \ldots , ( x_n, y_n ) $ \\

are repetition free in $X \times Y$. \\

{\bf Corollary 3.4.1.} \\

If an element $w \in X \bigotimes_{\lambda_X, \lambda_Y} Y$ is entangled, then it is of form (3.4.7), and the respective
repetition free (3.4.8) is of length at least two. \\

{\bf Remark 3.4.1.} \\

Obviously, the tensor products \\

$~~~~~~ X \bigotimes_{\lambda_X, \rho_Y} Y,~ X \bigotimes_{\rho_X, \lambda_Y} Y $ and $ X \bigotimes_{\rho_X, \rho_Y} Y $ \\

can be obtained in a similar manner. \\ \\

{\bf 3.5. A Fourth Example} \\

Let $X$ be a vector space over $\mathbb{R}$ or $\mathbb{C}$, and define on it the binary operation \\

(3.5.1)~~~ $ \mu ( x, x\,' ) = ( x + x\,' ) / 2,~~~ x, x\,' \in X $ \\

Similarly, on a vector space $Y$ over $\mathbb{R}$ or $\mathbb{C}$, we define the binary operation \\

(3.5.2)~~~ $ \nu ( y, y\,' ) = ( y + y\,' ) / 2,~~~ y, y\,' \in Y $ \\

Now we construct the tensor product \\

(3.5.3)~~~  $ X \bigotimes_{\mu, \nu} Y = Z / \approx_{\mu, \nu} $ \\

Particularizing the conditions (A2.1.6), (A2.1.7) which are involved in the definition of the equivalence relation $\approx_{\mu, \nu}$, we obtain that
(A2.1.6) takes the from \\

(3.5.4)~~~  replace $( x_1, y_1 ) ~\gamma~ ( x_2, y_2 ) ~\gamma~ \ldots ~\gamma~ ( x_n, y_n )$ \\
            \hspace*{1.85cm} with $( x\,'_1, y_1 ) ~\gamma~ ( x\,''_1, y_1 ) ~\gamma~ ( x_2, y_2 ) ~\gamma~ \ldots ~\gamma~ ( x_n, y_n )$, \\
            \hspace*{1.9cm}or vice-versa, where $x\,'_1, x\,''_1 \in X,~ x\,'_1 + x\,''_1 = 2 x_1$ \\

while (A2.1.7) becomes \\

(3.5.5)~~~  replace $( x_1, y_1 ) ~\gamma~ ( x_2, y_2 ) ~\gamma~ \ldots ~\gamma~ ( x_n, y_n )$ \\
            \hspace*{1.85cm} with $( x_1, y\,'_1 ) ~\gamma~ ( x_1, y\,''_1 ) ~\gamma~ ( x_2, y_2 ) ~\gamma~ \ldots ~\gamma~ ( x_n, y_n )$, \\
            \hspace*{1.9cm}or vice-versa, where $y\,'_1, y\,''_1 \in Y,~ y\,'_1 + y\,''_1 = 2 y_1$ \\

It follows that, for $x \in X,~ y \in Y$, we have \\

(3.5.6)~~~ $ \begin{array}{l}
                  x \bigotimes_{\mu, \nu} y ~=~ ( x\,' \bigotimes_{\mu, \nu} y ) ~\gamma~
                               ( x\,'' \bigotimes_{\mu, \nu} y ) ~=~ \\ \\
                  ~=~ ( x \bigotimes_{\mu, \nu} y\,' ) ~\gamma~ ( x \bigotimes_{\mu, \nu} y\,'' )
              \end{array} $ \\ \\

for $x\,', x\,'' \in X,~ x\,' + x\,'' = 2 x,~ y\,', y\,'' \in Y,~ y\,' + y\,'' = 2 y$. \\

Obviously, (3.5.6) is again a powerful simplification rule regarding the terms of the tensor product $X \bigotimes_{\mu, \nu} Y$. However, we can make
more clear the respective situation. Namely, any given finite sequence in $X \times Y$ \\

(3.5.7)~~~ $ ( x_1, y_1 ), ( x_2, y_2 ), ( x_3, y_3 ), \ldots , ( x_n, y_n ) $ \\

is called {\it median free}, if and only if \\

(3.5.8)~~~ $ 2 x_i \neq x_j + x_k,~~ 2 y_i \neq y_j + y_k $ \\

for $1 \leq i, j, k \leq n$, with $i, j, k$ pair-wise different. \\

Now, in view of (A2.1.10), we obtain \\

{\bf Theorem 3.5.1.} \\

The elements of the tensor product  $X \bigotimes_{\lambda_X, \lambda_Y} Y$ are of the form \\

(3.5.9)~~~ $ x_1 \bigotimes_{\lambda_X, \lambda_Y} y_1 ~\gamma~ x_2 \bigotimes_{\lambda_X, \lambda_Y} y_2 ~\gamma~
                                   \ldots ~\gamma~ x_n \bigotimes_{\lambda_X, \lambda_Y} y_n $ \\

where \\

(3.5.10)~~~ $ ( x_1, y_1 ), ( x_2, y_2 ), ( x_3, y_3 ), \ldots , ( x_n, y_n ) $ \\

are median free in $X \times Y$. \\

{\bf Corollary 3.5.1.} \\

If an element $w \in X \bigotimes_{\lambda_X, \lambda_Y} Y$ is entangled, then it is of form (3.5.9), and the respective repetition free (3.5.10) is of
length at least two. \\ \\

{\bf 4. Further Generalizations of Tensor Products} \\

The above example in subsection 3.5. gives an indication about both the possibility and interest in further generalizing tensor products, and doing so
beyond the already rather general setup in Appendix 2. \\

As a start, one can in a natural manner extend the concept of tensor products defined in (A2.1.8) which is based on binary operations. Instead, and as
suggested by (3.5.4), (3.5.5), one can use arbitrary multi-arity relations. Indeed, in (3.5.4), the $1 + 2$, or ternary relation on $X$, of the form \\

$~~~~~~ 2 x = x\,' + x\,'' $ \\

is used, and obviously this relation cannot be obtained as \\

$~~~~~~ x = \alpha ( x\,', x\,'' ) $ \\

from a binary operation $\alpha$ on $X$, in case $( X, \alpha )$ is, for instance, an arbitrary semigroup. Similar is the situation with (3.5.5). \\

Let therefore $X$ and $Y$ be two nonvoid sets, and let $A \subseteq X^n \times X^m$, while $B \subseteq Y^k \times Y^l$, for some $n, m, k, l \geq 1$.
Then, in extending (A2.1.5) - (A2.1.8), we define an equivalence relation $\approx_{A, B}$ on $Z$, as follows. Two sequences in $Z$ are equivalent, if
and only if they are identical, or each can be obtained from the other by a finite number of applications of the following operations \\

(4.1)~~~ permute pairs $( x_i, y_i )$ within the sequence \\

(4.2)~~~  replace $( x_1, y ) ~\gamma~ \ldots ~\gamma~ ( x_n, y ) ~\gamma~ ( u_1, v_1 ) ~\gamma~ \ldots ~\gamma~ ( u_p, v_p )$ \\
          \hspace*{1.85cm} with $( x\,'_1, y ) ~\gamma~ \ldots ~\gamma~ ( x\,'_m, y ) ~\gamma~ ( u_1, v_1 ) ~\gamma~ \ldots ~\gamma~ ( u_p, v_p )$
          \hspace*{1.9cm}or vice-versa, where $( x_1, \ldots , x_n, x\,'_1, \ldots , x\,'_m ) \in A$ \\

(4.3)~~~  replace $( x, y_1 ) ~\gamma~ \ldots ~\gamma~ ( x, y_k ) ~\gamma~ ( u_1, v_1 ) ~\gamma~ \ldots ~\gamma~ ( u_p, v_p )$ \\
          \hspace*{1.85cm} with $( x, y\,'_1 ) ~\gamma~ \ldots ~\gamma~ ( x, y\,'_l ) ~\gamma~( u_1, v_1 ) ~\gamma~ \ldots ~\gamma~ ( u_p, v_p )$, \\
          \hspace*{1.9cm}or vice-versa, where $( y_1, \ldots , y_k, y\,'_1, \ldots , y\,'_l ) \in B$ \\

Further details related to the resulting tensor products \\

(4.4)~~~ $ X \bigotimes_{A, B} Y = Z / \approx_{A, B} $ \\

will be presented elsewhere. \\

Here we note that in the example in subsection 3.4., we had \\

(4.5)~~~ $ A = \{~ ( x, x\,'_1, x\,'_2 ) \in X \times X^2 ~~|~~ 2 x = x\,'_1 + x\,'_2 ~\} $ \\

(4.6)~~~ $ B = \{~ ( y, y\,'_1, y\,'_2 ) \in Y \times Y^2 ~~|~~ 2 y = y\,'_1 + y\,'_2 ~\} $ \\

and obviously, we obtain $\approx_{A, B} ~=~ \approx_{\mu, \nu}$ on $Z$. \\

This example already shows the interest in the above generalization. Indeed, with $A, B$ as above, one can define the tensor product in subsection 3.4.
not only when $X$ and $Y$ are vector spaces, but also in the more general case when they are merely semigroups. And in such a case, the respective
equivalence $\approx_{A, B}$ is in general no longer of the form in (A2.2.1) - (A2.2.3), but of the more general form in (4.1) - (4.3) above. \\

Finally, the more general tensor products (A2.2.4) can also be further extended by considering families ${\cal A}$ of subsets $A \subseteq X^n \times
X^m$, respectively, families ${\cal B}$ of subsets $B \subseteq Y^k \times Y^l$, and accordingly, modifying the operations (A2.2.2), (A2.2.3). \\ \\

{\bf Appendix 1. ( Definition of Usual Tensor Products of \\
                   \hspace*{3cm} Vector Spaces)} \\

For convenience, we recall here certain main features of the usual tensor product of vector
spaces, and relate them to certain properties of Cartesian products. \\

Let $\mathbb{K}$ be a field and $E,F, G$ vector spaces over $\mathbb{K}$. \\

{\bf A1.1. Cartesian Product of Vector Spaces} \\

Then $E \times F$ is the vector space over $\mathbb{K}$ where the operations are given by \\

$~~~~~~ \lambda ( x, y ) + \mu ( u, v ) ~=~ ( \lambda x + \mu u, \lambda y + \mu v ) $ \\

for any $x, y \in E,~ u, v \in F,~ \lambda, \mu \in \mathbb{K}$. \\ \\

{\bf A1.2. Linear Mappings} \\

Let ${\cal L} ( E, F )$ be the set of all mappings \\

$~~~~~~ f : E ~\longrightarrow~ F $ \\

such that \\

$~~~~~~ f ( \lambda x + \mu u ) ~=~ \lambda f ( x ) + \mu f ( u ) $ \\

for $u, v \in E,~ \lambda, \mu \in \mathbb{K}$. \\ \\

{\bf A1.3. Bilinear Mappings} \\

Let ${\cal L} ( E, F; G )$ be the set of all mappings \\

$~~~~~~ g : E \times F ~\longrightarrow~ G $ \\

such that for $x \in E$ fixed, the mapping $F \ni y \longmapsto g ( x, y ) \in G$ is linear in
$y$, and similarly, for $y \in F$ fixed, the mapping $E \ni x \longmapsto g ( x, y ) \in G$ is
linear in $x \in E$. \\

It is easy to see that \\

$~~~~~~ {\cal L} ( E, F; G ) ~=~ {\cal L} ( E, {\cal L} ( F, G ) ) $ \\ \\

{\bf A1.4. Tensor Products} \\

The aim of the tensor product $E \bigotimes F$  is to establish a close connection between the
{\it bilinear} mappings in ${\cal L} ( E, F; G )$ and the {\it linear} mappings in
${\cal L} ( E \bigotimes F , G )$. \\

Namely, the {\it tensor product} $E \bigotimes F$ is : \\

(A1.4.1)~~~ a vector space over $\mathbb{K}$, together with \\

(A1.4.2)~~~ a bilinear mapping $t : E \times F ~\longrightarrow~ E \bigotimes F$, such that we \\
          \hspace*{1.6cm} have the following : \\ \\

{\bf UNIVERSALITY PROPERTY} \\

$\begin{array}{l}
    ~~~~~~~~~~ \forall~~~ V ~\mbox{vector space over}~ \mathbb{K},~~
                   g \in {\cal L} ( E, F; V ) ~\mbox{bilinear mapping}~ : \\ \\
    ~~~~~~~~~~ \exists~ !~~ h \in {\cal L} ( E \bigotimes F, V )
                                         ~\mbox{linear mapping}~ : \\ \\
    ~~~~~~~~~~~~~~~~ h \circ t ~=~ g
  \end{array} $ \\ \\

or in other words : \\

(A1.4.3)~~~ the diagram commutes

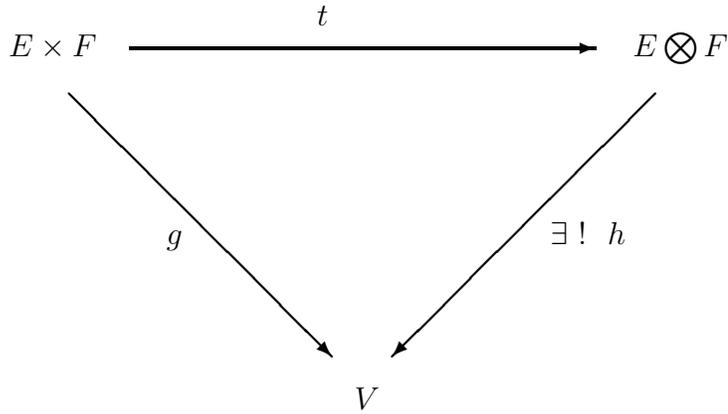
\begin{math}
\setlength{\unitlength}{1cm}
\thicklines
\begin{picture}(13,7)

\put(0.9,5){$E \times F$}
\put(2.5,5.1){\vector(1,0){6.2}}
\put(9.2,5){$E \bigotimes F$}
\put(5,5.4){$t$}
\put(1.7,4.5){\vector(1,-1){3.5}}
\put(9.5,4.5){\vector(-1,-1){3.5}}
\put(5.5,0.3){$V$}
\put(3,2.5){$g$}
\put(8.1,2.5){$\exists~!~~ h$}

\end{picture}
\end{math}

and \\

(A1.4.4)~~~ the tensor product $E \bigotimes F$ is {\it unique} up to vector \\
          \hspace*{1.6cm} space isomorphism. \\

Therefore we have the {\it injective} mapping \\

$~~~~~~ {\cal L} ( E, F; V ) \ni g ~\longmapsto~ h \in {\cal L} ( E \bigotimes F, V )
                                                    ~~~~\mbox{with}~~~~ h \circ t ~=~ g $ \\

The converse mapping \\

$~~~~~~ {\cal L} ( E \bigotimes F, V ) \ni h ~\longmapsto~
                         g ~=~ h \circ t \in {\cal L} ( E, F; V ) $ \\

obviously exists. Thus we have the {\it bijective} mapping \\

$~~~~~~ {\cal L} ( E \bigotimes F, V ) \ni h ~\longmapsto~
                         g ~=~ h \circ t \in {\cal L} ( E, F; V ) $ \\ \\

{\bf A1.5. Lack of Interest in ${\cal L} ( E \times F, G )$} \\

Let $f \in {\cal L} ( E \times F, G )$ and $( x, y ) \in E \times F$, then $( x, y ) = ( x, 0 )
+ ( 0, y )$, hence \\

$~~~~~~ f ( x, y ) ~=~ f ( ( x, 0 ) + ( 0, y ) ) ~=~ f ( x, 0 ) ~+~ f ( 0, y ) $ \\

thus $f ( x, y )$ depends on $x$ and $y$ in a {\it particular} manner, that is, separately on
$x$, and separately on $y$. \\ \\

{\bf A1.6. Universality Property of Cartesian Products} \\

Let $X, Y$ be two nonvoid sets. Their cartesian product is : \\

(A1.6.1)~~~ a set $X \times Y$, together with \\

(A1.6.2)~~~ two projection mappings~ $p_X : X \times X ~\longrightarrow~ X, \\
          \hspace*{1.6cm} p_Y : X \times Y ~\longrightarrow~ Y$, such that we have the
          following : \\ \\

{\bf UNIVERSALITY PROPERTY} \\

$\begin{array}{l}
    ~~~~~~~~~~ \forall~~~ Z ~\mbox{nonvoid set},~~ f : Z ~\longrightarrow~ X,~~
                          g : Z ~\longrightarrow~ Y~ : \\ \\
    ~~~~~~~~~~ \exists~ !~~ h  : Z ~\longrightarrow~ X \times Y~ : \\ \\
    ~~~~~~~~~~~~~~~~ f ~=~ p_X \circ h,~~~ g ~=~ p_Y \circ h
  \end{array} $ \\ \\

or in other words : \\

(A1.6.3)~~~ the diagram commutes

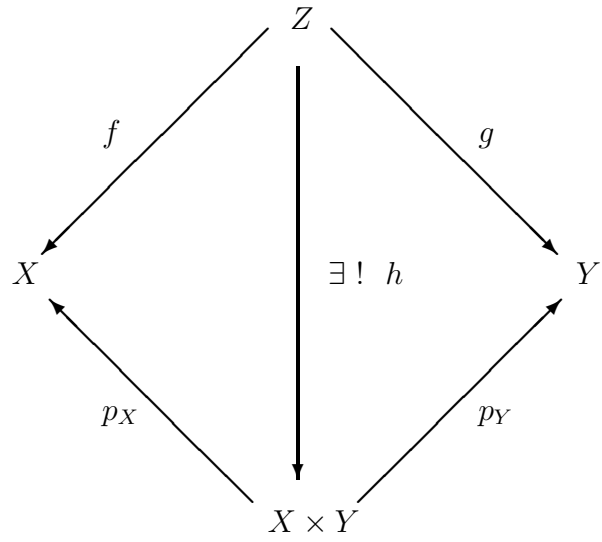
\begin{math}
\setlength{\unitlength}{1cm}
\thicklines
\begin{picture}(20,9)

\put(6.5,3.6){$\exists~ ! ~~ h$}
\put(6.1,6.5){\vector(0,-1){5.5}}
\put(6,7){$Z$}
\put(3.5,5.5){$f$}
\put(5.7,7){\vector(-1,-1){3}}
\put(8.5,5.5){$g$}
\put(6.55,7){\vector(1,-1){3}}
\put(2.3,3.6){$X$}
\put(9.8,3.6){$Y$}
\put(3.5,1.8){$p_X$}
\put(5.5,0.7){\vector(-1,1){2.7}}
\put(8.5,1.8){$p_Y$}
\put(6.9,0.7){\vector(1,1){2.7}}
\put(5.7,0.3){$X \times Y$}

\end{picture}
\end{math} \\ \\

{\bf A1.7. Cartesian and Tensor Products seen together} \\

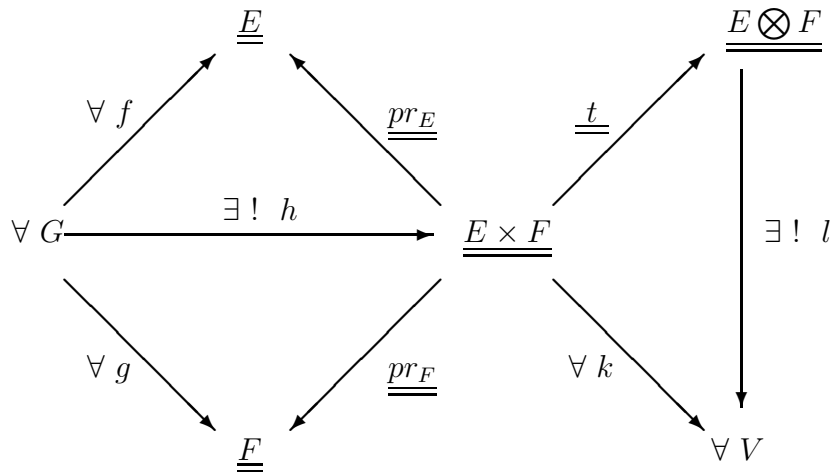
\begin{math}
\setlength{\unitlength}{1cm}
\thicklines
\begin{picture}(30,10)

\put(-0.2,6){$\forall~G$}
\put(0.5,6.5){\vector(1,1){2}}
\put(0.8,7.6){$\forall~f$}
\put(2.8,8.8){$\underline{\underline{E}}$}
\put(0.5,5.5){\vector(1,-1){2}}
\put(0.8,4.2){$\forall~g$}
\put(2.8,3.1){$\underline{\underline{F}}$}
\put(5.8,6){$\underline{\underline{E \times F}}$}
\put(5.5,6.5){\vector(-1,1){2}}
\put(4.8,7.6){$\underline{\underline{pr_E}}$}
\put(5.5,5.5){\vector(-1,-1){2}}
\put(4.8,4.2){$\underline{\underline{pr_F}}$}
\put(0.5,6.1){\vector(1,0){4.9}}
\put(2.6,6.3){$\exists~ !~~ h$}
\put(7,6.5){\vector(1,1){2}}
\put(7.3,7.6){$\underline{\underline{~t~}}$}
\put(7,5.5){\vector(1,-1){2}}
\put(7.2,4.2){$\forall~k$}
\put(9.3,8.8){$\underline{\underline{E \bigotimes F}}$}
\put(9.1,3.1){$\forall~V$}
\put(9.5,8.3){\vector(0,-1){4.5}}
\put(9.8,6){$\exists~ ! ~~ l$}

\end{picture}
\end{math} \\ \\

{\bf Appendix 2. ( Tensor Products beyond Vector Spaces)} \\

{\bf A2.1. The Case of Arbitrary Binary Operations} \\

Let us present the {\it first extension} of the standard definition of {\it tensor product} $X \bigotimes Y$, see Appendix 1, to the case of  two structures $( X, \alpha )$ and $( Y, \beta )$, where $\alpha : X \times X \longrightarrow X,~ \beta : Y \times Y \longrightarrow Y$ are arbitrary binary operations on two arbitrary given sets $X$ and $Y$, respectively. The way to proceed is as follows. Let us denote by $Z$ the set of all finite sequences of pairs \\

(A2.1.1)~~~ $ ( x_1, y_1 ), \dots , ( x_n, y_n ) $ \\

where $n \geq 1$, while $x_i \in X,~ y_i \in Y$, with $1 \leq i \leq n$. We define on $Z$ the binary operation $\gamma$ simply by the concatenation of the sequences (A2.1.1). It follows that $\gamma$ is associative, therefore, each  sequence (A2.1.1) can be written as \\

(A2.1.2)~~~ $ ( x_1, y_1 ), \dots , ( x_n, y_n ) =
                 ( x_1, y_1 ) ~\gamma~ ( x_2, y_2 ) ~\gamma~ \ldots ~\gamma~ ( x_n, y_n ) $ \\

where for $n = 1$, the right hand term is understood to be simply $( x_1, y_1 )$. Obviously, if $X$ or $Y$ have at least two elements, then $\gamma$ is not commutative. \\

Thus we have \\

(A2.1.3)~~~ $ Z = \left \{ ( x_1, y_1 ) ~\gamma~ ( x_2, y_2 ) ~\gamma~ \ldots ~\gamma~ ( x_n, y_n ) ~~
                             \begin{array}{|l}
                              ~ n \geq 1 \\ \\
                              ~ x_i \in X,~ y_i \in Y,~ 1 \leq i \leq n
                             \end{array} \right \} $ \\ \\

which clearly gives \\

(A2.1.4)~~~ $ X \times Y \subseteq Z $ \\

Now we define on $Z$ an equivalence relation $\approx_{\alpha, \beta}$ as follows. Two sequences in (A2.1.1) are equivalent, if and only if they are identical, or each can be obtained from the other by a finite number of applications of the following operations \\

(A2.1.5)~~~ permute pairs $( x_i, y_i )$ within the sequence \\

(A2.1.6)~~~ replace $( \alpha ( x_1, x\,'_1 ) , y_1 ) ~\gamma~ ( x_2, y_2 ) ~\gamma~ \ldots ~\gamma~ ( x_n, y_n )$ \\
            \hspace*{1.9cm} with $( x_1, y_1 ) ~\gamma~ ( x\,'_1, y_1 ) ~\gamma~ ( x_2, y_2 ) ~\gamma~ \ldots ~\gamma~
            ( x_n, y_n )$, \\
            \hspace*{1.9cm}or vice-versa \\

(A2.1.7)~~~ replace $( x_1, \beta ( y_1, y\,'_1 ) ) ~\gamma~ ( x_2, y_2 ) ~\gamma~ \ldots ~\gamma~ ( x_n, y_n )$ \\
            \hspace*{1.9cm} with $( x_1, y_1 ) ~\gamma~ ( x_1, y\,'_1 ) ~\gamma~ ( x_2, y_2 ) ~\gamma~ \ldots ~\gamma~
            ( x_n, y_n )$, \\
            \hspace*{1.9cm}or vice-versa \\

Finally, the {\it tensor product} of $( X, \alpha )$ and $( Y, \beta )$ is defined to be the
quotient space \\

(A2.1.8)~~~ $ X \bigotimes_{\alpha, \beta} Y = Z / \approx_{\alpha, \beta} $ \\

with the canonical quotient embedding, see (A2.1.4) \\

(A2.1.9)~~~ $ X \times Y \ni ( x, y ) \longmapsto x \bigotimes_{\alpha, \beta} y \in X \bigotimes_{\alpha, \beta} Y $ \\

where as in the usual case of tensor products, we denote by $x \bigotimes_{\alpha, \beta} y$, or simply  $x \bigotimes y$,
the equivalence class of $( x, y ) \in X \times Y \subseteq Z$. \\

Obviously, the binary operation $\gamma$ on $Z$ will canonically lead by this quotient operation to a {\it commutative}
and {\it associative} binary operation on $X \bigotimes_{\alpha, \beta} Y$, which for convenience is denoted by the same
$\gamma$, although in view of (A2.1.8), this time it depends on $\alpha$ and $\beta$, thus it should rigorously be written
$\gamma_{\alpha, \beta}$. \\

In this way, the elements of $X \bigotimes_{\alpha, \beta} Y$ are all the expressions \\

(A2.1.10)~~~ $ x_1 \bigotimes_{\alpha, \beta} y_1 ~\gamma~ x_2 \bigotimes_{\alpha, \beta} y_2 ~\gamma~
                                   \ldots ~\gamma~ x_n \bigotimes_{\alpha, \beta} y_n $ \\

with $n \geq 1$ and $x_i \in X,~ y_i \in Y$, for $1 \leq i \leq n$. \\

The customary particular situation is when $X$ and $Y$ are commutative semigroups, groups, or even vector spaces over some field $\mathbb{K}$. In this case $\alpha, \beta$ and $\gamma$ are as usual denoted by +, that is, the sign of addition. \\

It is easy to note that in the construction of tensor products above, it is {\it not} necessary for $( X, \alpha )$ or
$( Y, \beta )$ to be semigroups, let alone groups, or for that matter, vector spaces. Indeed, it is sufficient that
$\alpha$ and $\beta$ are arbitrary binary operations on $X$ and $Y$, respectively, while $X$ and $Y$ can be arbitrary
sets. \\

Also, as seen above, $\alpha$ or $\beta$ need {\it not} be commutative either. However, the resulting tensor product
$X \bigotimes_{\alpha, \beta} Y$, with the respective binary operation $\gamma$, will nevertheless be commutative
and associative. \\ \\

{\bf A2.2. A Second Generalization of Tensor Products} \\

The first generalization of tensor products presented in section A2.1. above can easily be further extended. Indeed, let
$X,~ Y$ be arbitrary sets and let ${\cal A}$ be any set of binary operations on $X$, while correspondingly, ${\cal B}$ is
any set of binary operations on $Y$. \\

The constructions in (A2.2.1) - (A2.1.4) can again be implemented, since they only depend on the sets $X,~ Y$. \\

Now, we can define on $Z$ the equivalence relation $\approx_{{\cal A}, {\cal B}}$ as follows. Two sequences in (A2.1.1)
are equivalent, if and only if they are identical, or each can be obtained from the other by a finite number of
applications of the following operations \\

(A2.2.1)~~~ permute pairs $( x_i, y_i )$ within the sequence \\ \\

(A2.2.2)~~~ replace $( x_1, y_1 ) ~\gamma~ ( x\,'_1, y_1 ) ~\gamma~ ( x_2, y_2 ) ~\gamma~ \ldots ~\gamma~ ( x_n, y_n )$ \\
            \hspace*{1.9cm} with $( \alpha ( x_1, x\,'_1), y_1 ) ~\gamma~ ( x_2, y_2 ) ~\gamma~ \ldots , ( x_n, y_n )$, \\
            \hspace*{1.9cm} or vice-versa, where $\alpha \in {\cal A}$ \\ \\

(A2.2.3)~~~ replace $( x_1, y_1 ) ~\gamma~ ( x_1, y\,'_1 ) ~\gamma~ ( x_2, y_2 ) ~\gamma~ \ldots ~\gamma~ ( x_n, y_n )$ \\
            \hspace*{1.9cm} with $( x_1, \beta ( y_1, y\,'_1 ) ) ~\gamma~ ( x_2, y_2 ) ~\gamma~ \ldots ~\gamma~
            ( x_n, y_n )$, \\
            \hspace*{1.9cm} or vice-versa, where $\beta \in {\cal B}$ \\ \\

Thus one obtains the tensor product \\

(A2.2.4)~~~ $ X \bigotimes_{{\cal A}, {\cal B}} Y = Z / \approx_{{\cal A}, {\cal B}} $ \\

with the canonical quotient embedding \\

(A2.2.5)~~~ $ X \times Y \ni ( x, y ) \longmapsto x \bigotimes_{{\cal A}, {\cal B}} y \in X \bigotimes_{{\cal A},
{\cal B}} Y $ \\

and $X \bigotimes_{{\cal A}, {\cal B}} Y$ being the set of all elements \\

(A2.2.6)~~~ $ x_1 \bigotimes_{{\cal A}, {\cal B}} y_1 ~\gamma~ x_2 \bigotimes_{{\cal A}, {\cal B}} y_2 ~\gamma~
                                  \ldots ~\gamma~ x_n \bigotimes_{{\cal A}, {\cal B}} y_n $ \\

with $n \geq 1$ and $x_i \in X,~ y_i \in Y$, for $1 \leq i \leq n$. \\

Let us further note that, given sets ${\cal A} \subseteq {\cal A}\,'$ of binary operations on $X$, and correspondingly,
sets ${\cal B} \subseteq {\cal B}\,'$ of binary operations on $Y$, we have the surjective mapping \\

(A2.2.7)~~~ $ X \bigotimes_{{\cal A}, {\cal B}} Y \ni ( z )_{{\cal A}, {\cal B}} \longmapsto
                         ( z )_{{\cal A}\,', {\cal B}\,'} \in  X \bigotimes_{{\cal A}\,', {\cal B}\,'} Y $ \\

where $( z )_{{\cal A}, {\cal B}}$ denotes the $\approx_{{\cal A}, {\cal B}}$ equivalence class of $z \in Z$, and
similarly with $( z )_{{\cal A}\,', {\cal B}\,'}$. \\ \\

\end{document}